\documentclass[11pt]{article}
\usepackage[utf8]{inputenc}
\usepackage[T1]{fontenc}
\usepackage{amssymb}
\usepackage{amsfonts}
\usepackage{floatrow}
\usepackage{amsthm}
\usepackage{graphicx}
\usepackage{amsmath}
\usepackage{eucal}
\usepackage{stmaryrd}
\usepackage{caption}
\usepackage{mathtools}
\usepackage{amsfonts,epsfig,epstopdf,titling,url,array}

\numberwithin{equation}{section}

\newtheorem{defn}{Definition}[section]
\newtheorem{conj}{Conjecture}[section]

\newtheorem{thm}{Theorem}
\newtheorem{prop}{Proposition}
\newtheorem{lemma}{Lemma}
\newtheorem{corr}{Corollary}
\newcommand{\beq}{\begin{equation}}
\newcommand{\RN}[1]{%
  \textup{\uppercase\expandafter{\romannumeral#1}}%
}
\newcommand{\eeq}{\end{equation}}
\newcommand{\integral}{\int_{\partial\Omega}}
\newcommand{\nablat}{\nabla_{\tau}}
\newcommand{\beqs}{\begin{equation*}}
\newcommand{\eeqs}{\end{equation*}}

\usepackage[
top    = 2.75cm,
bottom = 2.50cm,
left   = 2.3cm,
right  = 2.3cm]{geometry}

\title{A Critical Domain For the First Normalized Nontrivial Steklov Eigenvalue Among Planar Annular Domains} 

\author{
\small Leoncio Rodriguez Quiñones, email: rodriguez-l@javeriana.edu.co\\ \small Departamento de Matem\'{a}ticas\\
\small Pontificia Universidad Javeriana, Bogot\'{a}, Colombia
}

\begin{document}
\bibliographystyle{plain}
\large{





\maketitle

\begin{abstract}
We describe a shape derivative approach to provide a candidate for an optimal domain among non-simply connected planar domains with two boundary components. This approach is an adaptation of the work on the extremal eigenvalue problem for the Wentzell-Laplace operator developed by Dambrine, Kateb and Lamboley \cite{Dambrine}.
\end{abstract}

\section{Introduction}

The problem of finding domains that maximize or minimize a given eigenvalue (or functions of eigenvalues) of an elliptic operator turns out to be difficult to solve. The challenge usually comes from the fact that to prove statements about their existence, or regularity properties, or even to study their shapes; it is necessary to borrow tools from different areas of mathematics. In particular, areas such as calculus of variations, differential geometry, analysis and partial differential equations turn out to be useful when working on this type of problems.


Even for classical eigenvalue problems such as the Dirichlet eigenvalue problem and Neumann eigenvalue problem, there are still a several open problems related to the spectral geometry for lower eigenvalues (see \cite{AHen06}).  

Historically, it is well-known that one of the first works on shape optimization of eigenvalues of boundary value problems was due to  Lord Rayleigh. Rayleigh conjectured that for simply connected bounded planar domains, the first Dirichlet eigenvalue of the Laplacian satisfies the inequality 
$\lambda_{1}(\Omega)|\Omega|^{1/2}\geq \sqrt[]{\pi}\lambda_{1}(\mathbb{D})$. Rayleigh presented several cases where equality is achieved when $\Omega=\mathbb{D}$ and  $\lambda_{1}(\mathbb{D})$ is the least positive zero of the Bessel function $J_{0}(r)$ \cite{rayleigh}. In $1923$ and $1924$ Faber and Krahn proved independently Rayleigh's conjecture, their result is now known as the Faber-Krahn inequality (\cite{faber}, \cite{krahn}). Almost thirty years later, in $1951$, P\'{o}lya and Sz\"{e}go also showed that the inequality $\lambda_{1}(\Omega)|\Omega|^{1/2}\geq \sqrt[]{\pi}\lambda_{1}(\mathbb{D})$ can be refined by using Steiner symmetritation on $\Omega$ (\cite{polya}, \cite{szego}).

In the case of the Neumann eigenvalues, Kornhauser and Stakgold in $1952$, conjectured  that among domains $\Omega$ of prescribed area, the circle maximizes $\mu_{1}(\Omega)|\Omega|\leq\mu_{1}(\mathbb{D})\pi$, where $\mu_{1}$ is the first simple nontrivial eigenvalue \cite{stakgold}. It was later in $1954$, that Sz\"{e}go proved Kornhauser and Stakgold's conjecture in \cite{szego}, and in $1956$ Weinberger generalized Sz\"{e}go's result to the $N-$dimensional case \cite{weinberger}.

It is worth metioning that in the case  of the Neumann eigenvalues, the minimization problem happens to be easy and it follows  from  the construction of a rather simple domain, namely a rectangle. Thus, for Neumann eigenvalues the interesting problems are about the domains that maximize them (see \cite{AHen06}, Ch. 7.). 

In this note, our focus is on a different boundary condition for the Laplacian. We are interested in the so called \textit{Steklov} eigenvalue problem, given by 
\begin{align}
     -\Delta u &=0,\quad x\in\Omega,\\ \nonumber
      \partial_{n}u &= \lambda u,\quad  x\in\partial\Omega, \label{equ1}
\end{align} 
where $\partial_{n}u$ is the outer normal derivative and $u\not\equiv 0$. The Steklov problem arises in the modeling of the vibration of a free membrane whose whole mass is uniformly distributed on $\partial\Omega$  \cite{Girouard1}.


The Steklov eigenvalues satisfy the following variational characterization,
\begin{equation}\label{equ2}
\lambda_{n}(\Omega)=\min_{u\in H^{1}(\Omega)}\left\{\frac{\int_{\Omega}|\nabla u|^{2}dx}{\int_{\partial\Omega}u^{2}ds}:\int_{\partial\Omega} u\phi_{j}ds=0, j=0,1,\dots,n-1\right\},
\end{equation} 
where $\phi_{j}$ is the eigenfunction associated to the $j$-th eigenvalue (see for example, \cite{bandle},\cite{Brock},\cite{AHen06},\cite{Girouard1}). Since the Steklov eigenvalues correspond to the eigenvalues of a positive, formally self-adjoint pseudodifferential operator of order one, its spectrum satisfies
\begin{align*}
0=\lambda_{0}(\Omega)<\lambda_{1}(\Omega)\leq\lambda_{2}(\Omega)\leq \dots \leq\lambda_{n}(\Omega)\leq\lambda_{n+1}(\Omega)\leq\dots\rightarrow\infty
\end{align*}
where the eigenvalue $\lambda_{0}(\Omega)$ corresponds to the constant eigenfunction \cite{Dambrine}. 

As in the Neumann eigenvalues case, the  maximization problem for the Steklov eingenvalues is the more interesting one. Although the construction of a minimizing domain is not as easy as in the Neumann case, it is still possible to find one for which $\lim_{\epsilon\rightarrow 0}\lambda_{n}(\Omega_{\epsilon})=0$. In this case,  $\Omega_{\epsilon}$ is constructed by joining two copies of the unit disk with a rectangle of width $\epsilon^{3}$, and length $\epsilon$ \cite{Girouard4}.

Thus, when considering the maximization problem for the Steklov eigenvalues, a well-known result for planar domains is given by Weinstock's inequality, which states that in the class of simply-connected bounded domains, the unit disk $\mathbb{D}$ is a maximizer for the normalized first Steklov eigenvalue $\lambda_{1}(\Omega)|\partial\Omega|$ \cite{weinstock}. This normalization is considered, so the problem of maximizing this product is equivalent to maximize $\lambda_{n}(\Omega)$ with a fixed perimeter constraint on $\Omega$ \cite{Girouard3}. The generalization of Weinstock's result to higher dimensions was proved by \textit{Brock} \cite{Brock}. It states that the ball maximizes the first nontrivial Steklov eigenvalue among open sets of given volume.

It is important to mention that the assumption of simply-connectivity on the domain needs to be kept in order for Weinstock's result to hold. One example that shows that Weinstock's inequality fails if such assumption is removed, is given by annuli. 

The Steklov eigenvalues on an annulus $\Omega_{\epsilon}$ with outer radius $r_{o}=1$ and inner radius $r_{i}=\epsilon$ are given by
\begin{align}\label{steklovevals}
\lambda_{n}(\Omega_{\epsilon})=\frac{n}{2}\left(\frac{1+\epsilon}{\epsilon}\right)\left(\frac{1+\epsilon^{2n}}{1-\epsilon^{2n}}\right)\pm\frac{n}{2}\sqrt{\left(\frac{1+\epsilon}{\epsilon}\right)^{2}\left(\frac{1+\epsilon^{2n}}{1-\epsilon^{2n}}\right)^{2}-\frac{4}{\epsilon}}.
\end{align}
 Taking 
$\epsilon$ small enough in (\ref{steklovevals}), it follows that  $\lambda_{1}(\Omega_{\epsilon})|\partial\Omega_{\epsilon}|>2\pi\sigma_{1}(\mathbb{D})$, showing that Weinstock's inequality fails \cite{Girouard3}. 

Knowing that Weinstock's inequality does not hold on bounded planar annular domains; one may ask if among such domains, there exists a domain  (or domains) that provides a maximizer for the first normalized Steklov eigenvalue. The conjecture is that among bounded planar domains with one hole (that is, domains whose complement has two connected components one of which is bounded and the other one unbounded) and Lipschitz boundary, the annulus whose inner radius determines the maximum of the function $\lambda_{1}(\Omega_{\epsilon})|\partial\Omega_{\epsilon}|$ is the maximizer. In other words the annulus whose inner radius gives the maximum of the curve showed in Figure \ref{normalizedfirst}.
\begin{figure}[ht]
\includegraphics[scale=0.6]{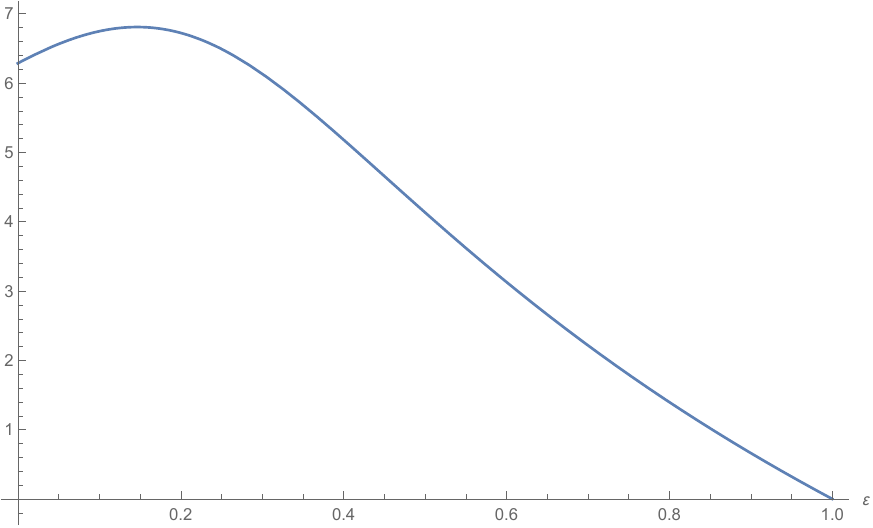}
\centering
\caption{Normalized eigenvalue $\sigma_{1}(\Omega_{\epsilon})|\partial\Omega_{\epsilon}|$. The max is attained at the solution $\epsilon_{0}$ of $\frac{d(\sigma_{1}(\Omega_{\epsilon})|\partial\Omega_{\epsilon}|)}{d\epsilon}=0.$ Numerically we get $\epsilon_{0}\approx0.146721$}
\label{normalizedfirst}
\end{figure}
This note shows that in a local sense, the annulus with outer radius $r_{o}=1$ and inner radius specified above, provides a critical domain for $\lambda_{1}(\Omega)|\partial\Omega|$ in the class of bounded planar domains with one hole, fixed outer boundary given by $\mathbb{S}^{1}$ and sufficiently smooth inner boundary. As a consequence we have that such critical domain is a candidate for a maximizer of such shape functional. The fact that this particular annulus is a critical domain, suggests that the candidate for maximizer could be symmetric, and with well known topological and geometrical properties.

This paper is structured as follows: section 2 presents the main results used to prove our assertion. Section 3 is focused on applying the theory presented in section 2 to our particular problem in order to prove Theorem 5, which is the main result of this paper, the computations presented in this section are an adaption to those presented in \cite{Dambrine}. We also want to point out that the proofs of Lemma 2, Propositions 1, 2 and 3 in \cite{gomes} were useful at the moment of deriving and understanding some these computations \cite{gomes}. Finally, in section 4 we present some numerical examples to verify the main result.
\section{Shape Calculus}\label{secshapecal}

In this section we present a brief introduction to the theory of shape calculus. We mainly focus on the tools needed to prove that an annulus provides a critical domain for the first nontrivial Steklov eigenvalue. All of the definitions and theorems needed to prove such statement, as well as a detailed presentation about the theory of shape calculus and geometries can be found in \cite{delfour}, \cite{sokolowski}. More general variation formulas on smooth manifolds for the eigenvalues of the Laplace-Beltrami operator can be found in \cite{berger}. For an interesting application of Hadamard type variation formulas for the eigenvalues of the so called $\eta$-Laplacian we refer the reader to \cite{gomes}. And for generic properties of eigenvalues and eigenfunctions, Uhlenbeck's \cite{uhlenbeck} work is of upmost importance.

We begin our presentation by introducing the concept of shape functional which is the function that we want to differentiate in some sense.
 
\begin{defn}
Given a nonempty subset $D$ of $\mathbb{R}^N$ and consider its power set $\CMcal{P}(D)$. A \textit{shape functional} is a map 
\begin{align*}
J:\CMcal{A}\rightarrow \mathbb{R}
\end{align*}
from some \textit{admissible family} $\CMcal{A}$ of sets in $\CMcal{P}(D)$ into $\mathbb{R}$.  The set $D$ will be referred to as the underlying \textit{holdall}.
\end{defn}
The simplest examples of shape functionals that we encounter are given by 
\begin{align*}
J(\Omega)=|\Omega| &= \int_{\Omega}dx, \\
J(\Omega)=|\partial\Omega| &= \int_{\partial\Omega}dS.
\end{align*}

In order to study variations of a shape functional, we want to consider perturbations of the  given domain $\Omega$. This perturbation will be defined by means of the transformation $T_{t}:\overline{\Omega}\rightarrow \overline{\Omega}_{t}$ defined by 
$$T_{t}(x)\coloneqq x+tV(x),$$
where the vector field satisfies $V\in W^{3,\infty}(\Omega,\mathbb{R}^{N})$ (see, Theorem 3.4 in \cite{Dambrine}, Prop. 3 in \cite{gomes}).


The condition on $V$  ensures that the one parameter family of domains $(\Omega_{t})_{t\geq 0}$ remains in the admissible family $\CMcal{A}$ that is being considered. In particular, these vector fields are meant to preserve the topological assumptions made on the original domain $\Omega$.

Once we have established the perturbations of the original domain $\Omega$, we need to introduce the concept of Eulerian derivative. This derivative will allow us to study how the given shape functional varies through the family of perturbed domains. 
\begin{defn}\label{defeulerder}
For any vector field $V \in W^{3,\infty}(\Omega,\mathbb{R}^{N})$ , the Eulerian derivative of the domain functional $J(\Omega)$ at $\Omega$ in the direction of a vector field $V$ is defined as the limit 
\begin{equation}
dJ(\Omega;V)=\lim_{t\rightarrow 0}\frac{J(\Omega_{t})-J(\Omega)}{t}
\end{equation}
where $\Omega_{t}=T_{t}(V)(\Omega)$.
\end{defn}

From this definition it is possible to conclude that for a function $u$ lying in a suitable function space, the functional
\begin{equation}
J(\Omega_{t})=\int_{\partial\Omega}u(t,x)dS
\end{equation}
has Eulerian derivative given by 
\begin{equation}\label{shapederfor}
dJ(\Omega;V) = \int_{\partial\Omega}u^{\prime}(\Omega;V)dS+\int_{\partial\Omega}\left(\frac{\partial u}{\partial n}+Hu\right)V_{n}dS,
\end{equation}
where $u^{\prime}=\partial_{t}u(t,x)|_{t=0}$ is the shape derivative, $V_{n}$ is the normal component of the vector field $V$ at $t=0$ , $H=\Delta b$ is the mean curvature on $\partial\Omega$ and $b$ is the oriented distance to $\Omega$.

It is important to mention that we are considering the derivative of the shape functional in only one direction, namely, the direction given by the vector field $V$.

\subsection{Shape derivative of eigenvalues and eigenfunctions}\label{section4}

To study how the eigenvalues and eigenfunctions change with respect to changes in the original domain, we need to study first regularity properties of the related quantities. We use most of the results related to the Wentzell boundary value problem given by

\begin{align}\label{wentzello}
     -\Delta u &=0,\quad x\in\Omega,\\ \nonumber
      -\beta\Delta_{\tau}u+\partial_{n}u &= \lambda u, \quad x\in\partial\Omega.
\end{align} 

A completed study and proofs of statements related with the Wentzell eigenvalues and eigenfunctions is presented in \cite{Dambrine}. 

For the sake of completeness, we present the most relevant results regarding the regularity of the eigenvalues and eigenfunctions of the Wentzell eigenvalue problem. The following Theorem guarantees that even though a multiple eigenvalue of an elliptic operator is not shape differentiable (see, \cite{sokolowski}), one still can compute the shape derivatives of their ``branches''. These branches arise as a result of perturbations of the domain by the vector field $V$. 

\begin{thm}[Dambrine-Kateb-Lamboley, \cite{Dambrine}] \label{thm1}
Let $\Omega$ be an open smooth bounded domain of $\mathbb{R}^{N}$, $V\in W^{3,\infty}(\Omega,\mathbb{R}^{N})$, $\Omega_{t}=T_{t}(\Omega)$, and $\lambda$ a multiple eigenvalue of the problem
\begin{align}\label{wentzell}
     -\Delta u &=0,\quad x\in\Omega,\\ \nonumber
      -\beta\Delta_{\tau}u+\partial_{n}u &= \lambda u, \quad x\in\partial\Omega.\label{equ1}
\end{align} 
Then, there exists $m$ (the multiplicity of the eigenvalue) real-valued continuous functions $t \rightarrow \lambda_{i}(t)$ and $m$ functions $t\rightarrow u^{t}_{i}\in H^{5/2}(\Omega)$, that are analytic in a neighborhood of $t=0$, $\lambda_{i}(0)=\lambda$ for $i = 1, \dots, m$, where the functions $u^{t}_{i}$ are the normalized eigenfunctions associated to the eigenvalue $\lambda_{i}(t)$ on the domain $\Omega_{t}$.
\end{thm}

Theorem \ref{thm1} therefore guarantees the needed regularity properties of eigenfunctions and eigenvalues, which allows us to consider their shape derivatives. As the next Theorem shows, it is possible to derive a boundary value problem whose solution is determined by the shape derivative of the eigenfunction solving problem \eqref{wentzell}. 

\begin{thm}[Dambrine-Kateb-Lamboley, \cite{Dambrine}] \label{thm2}
For an eigenpair path $(\lambda(t), u_{t})$ of $\Omega_{t}$, the shape derivative $u^{\prime}=(\partial_{t=0}u_{t})|_{t=0}$ of the eigenfunction $u_{t}$ for the Wentzell problem satisfies
\begin{align}
	\Delta u^{\prime} &=0, \quad x\in\Omega, \nonumber\\ 
	-\beta\Delta_{\tau}u^{\prime}+\partial_{n}u^{\prime}-\lambda u^{\prime} 	&= \beta\Delta_{\tau}(V_{n}\partial_{n}u)-\beta \text{div}_{\tau}(V_{n}			(2D^{2}b-HId)\nabla_{\tau}u) \\ 
	&+ \text{div}_{\tau}(V_{n}\nabla_{\tau}u)+\lambda^{\prime}(0)u+\lambda V_{n}(\partial_{n}u+Hu), \quad x\in\partial\Omega. \nonumber
\end{align}
\end{thm}

Notice that the boundary value problem in Theorem \ref{thm2} involves the shape derivative of the eigenvalue. From this expression then it is possible to obtain explicit formulas for $\lambda^{\prime}(0)$ for both cases when $\lambda$ is a simple or a multiple eigenvalue. 

For the case when $\lambda$ is a simple eigenvalue, multiplying both sides of the boundary condition in Theorem \ref{thm2} by the normalized eigenfunction $u$ and taking integrals in both sides of the equality gives the proof to the following Theorem.

\begin{thm}[Dambrine-Kateb-Lamboley, \cite{Dambrine}]\label{simplecase}
If $\lambda$ is a simple eigenvalue and $u$ the corresponding normalized eigenfunction then, the map $t\rightarrow \lambda(t)$ is analytic and its derivative at $t=0$ is
\begin{align*}
\lambda^{\prime}(0)=\int_{\partial\Omega}V_{n}\left(|\nabla_{\tau}u|^{2}-|\partial_{n}u|^{2}-\lambda H|u|^{2}+\beta (HId-2D^{2}b)\nabla_{\tau}u\cdot\nabla_{\tau}u\right)d\sigma.
\end{align*}
\end{thm}

As mentioned before, when the eigenvalue $\lambda$ is of multiplicity $m$, we no longer have differentiability. However, from Theorem \ref{thm1}, we know that for $t$ small enough, on the perturbed domain $\Omega_{t}$ there will be $m$ eigenvalues $\lambda_{i}(t)$, $i = 1,\dots, m$ for which it is possible to compute the derivative in the sense given by the Theorem \ref{multievals}.

\begin{thm}[Dambrine-Kateb-Lamboley, \cite{Dambrine}]\label{multievals}
Let $\lambda$ be a multiple eigenvalue of order $m \geq 2$. Then each $t \mapsto \lambda_{i}(t)$ for $i \in \llbracket 1,N \rrbracket $ given by Theorem \ref{thm1} has a derivative near $0$, and the values of $(\lambda_{i}^{\prime}(0))_{i \in \llbracket 1,N \rrbracket }$ are the eigenvalues of the matrix $M(V_{n})=(M_{jk})_{1 \leq j,k \leq m}$ defined by
\begin{align*}
M_{jk}=\int_{\partial \Omega}V_{n} \left( \nabla_{\tau}u_{j} \cdot \nabla_{\tau}u_{k}-\partial_{n}u_{j}\partial_{n}u_{k}-\lambda Hu_{j}u_{k}+\beta\left(HI_{d}-2D^{2}b\right) \nabla_{\tau}u_{j} \cdot \nabla_{\tau}u_{k}\right)d\sigma.
\end{align*}
\end{thm}
 
The proof of Theorem \ref{multievals} is similar to the proof of Theorem \ref{simplecase}. The important component of the proof is to notice that the solution to the problem \eqref{equ1} on $\Omega_{t}$ can be decomposed as $u=\sum_{j=1}^{m}c_{j}u_{j}$ where each $u_{j}$ are the associated eigenfunctions in the associated eigenspace.

Theorem \ref{multievals} also implies that balls in $\mathbb{R}^{N}$ are critical domains for the Wentzell eigenvalue. The following corollary present an explicit formula for the entries of the matrix $M$.

\begin{corr}[Dambrine-Kateb-Lamboley, \cite{Dambrine}]\label{balls}
If $\Omega$ is a ball of radius $R$ and $\lambda$ is the first non-trivial Wentzell eigenvalue, with multiplicity $n$. The shape derivative of the maps $t \mapsto \lambda_{i}(t)$, $i = 1,\dotso,n$ given by Theorem \ref{thm1} are the eigenvalues of the matrix $M_{B_{R}}(V_{n})=(M_{jk})_{j,k=1,\dotso,n}$ defined by
\begin{equation}\label{ballmatrix}
M_{jk}=\frac{\delta_{jk}}{\omega_{n}R^{n+1}}\left( 1+\beta\frac{n-3}{R}\right) \int_{\partial B_{R}}V_{n}-C(n,R)\int_{\partial B_{R}}V_{n}x_{j}x_{k}d\sigma,
\end{equation}
where $C(n,R)=\frac{(n+1)(1+\beta\frac{n-2}{R})}{\omega_{n}R^{n+3}}$.
\end{corr}
 
Knowing that the eigenvalues of the matrix defined by equation \eqref{ballmatrix} are the shape derivatives of the branches of the first Wentzell eigenvalue; the next step is to identify under what conditions these eigenvalues are either all zero, or they satisfy certain condition under which it is possible to conclude that a ball provide a critical shape. 

A vector field $V$ is said to be volume preserving at first order if $\small{\int_{\partial\Omega}V_{n}d\sigma=0}$. i.e if the normal component of the vector field is orthogonal to constants in $L^{2}(\partial\Omega)$. The next proposition states that if the deformation of the domain is produced by a volume preserving vector field then the eigenvalues of the matrix in \eqref{ballmatrix} are all zero. 

\begin{prop}[Dambrine-Kateb-Lamboley, \cite{Dambrine}]\label{zeromatrix}
If $V$ preserves volume, then the following statements are equivalent:
\begin{enumerate}
\item$V_{n}$ is orthogonal (in $L^{2}(\partial B(0,R))$) to homogeneous harmonic polynomials of degree 2, 
\item $M_{B(0,R)}(V_{n})=0$.
\end{enumerate}
\end{prop}

If the matrix defined by formula \eqref{ballmatrix} is not the null matrix, it is necessary to introduce the subgradient $\partial\lambda$ of $\lambda$, given by $\partial\lambda=[\inf_{i=0,\dots,n}\lambda^{\prime}_{i}(0),\sup_{i=0,\dots,n}\lambda^{\prime}_{i}(0)]$, to conclude that under volume preserving deformations; balls are still critical domains. The following proposition makes clear the need of introducing this subgradient.

\begin{prop}[Dambrine-Kateb-Lamboley, \cite{Dambrine}]\label{notzero}
When $\Omega$ is a ball of radius $R$, then
\begin{equation}
Tr(M_{B(0,R)}(V_{n}))=0
\end{equation}
for all volume preserving deformations.
\end{prop}

Thus, if the vector field is volume preserving but the matrix is not null, Proposition \ref{notzero} states that at least the sum of the eigenvalues is zero. This implies that at least one of the branches of the multiple eigenvalue is increasing and one is decreasing, so the ball is a critical domain for the maximum of the eigenvalue branches.

\section{The normalized eigenvalue}
We are now ready to show, that among the class of Lipschitz bounded planar domains with one hole, an annulus is a critical domain for the normalized first nontrivial Steklov eigenvlue.

To prove this, we take advantage of the results presented in section \ref{section4} regarding the first Wentzell eigenvalue with $\beta\in\mathbb{R}$. Once these results are established we obtain our results by considering the particular case $\beta=0$.

\subsection{Derivatives of normalized eigenvalues}
We begin by taking $\CMcal{A}$ as the class  of bounded, connected, Lipschitz, planar domains with genus zero and two boundary components. Intuitively, $\CMcal{A}$ is the class of bounded planar domains with one hole. We also consider the shape functional $J:\CMcal{A}\rightarrow [0,\infty]$ defined by $J(\Omega)=|\partial\Omega|\lambda(\Omega)$, where $|\cdot|$  denotes the Lebesgue measure of $\partial\Omega$.

Our first goal is to obtain an expression for the Eulerian derivative of $J$. Note that the functional $J$ is the product of two shape functionals, both of which have derivative in the Eulerian sense (for the derivative of $|\cdot|$, see \cite{sokolowski}, page 54). Our first step is to prove that for the functional $J$, the product rule holds in the sense of Eulerian differentiation. That the functional $J$ defined above has Eulerian derivative is a consequence of the  following lemma.

\begin{lemma}\label{productrule}
Let $I$ and $H$ be two shape functionals defined on a class $\CMcal{A}$ of admissible domains, both having Eulerian derivatives $dI(\Omega,V)$ and $dH(\Omega,V)$ respectively at $\Omega$, in the direction of a vector field  $V \in  W^{3,\infty}(\Omega,\mathbb{R}^{N})$. Then the shape functional $J$, defined by $J(\Omega)=I(\Omega)H(\Omega)$, has Eulerian derivative at  $\Omega$ in the direction of the vector field  $V$ and $dJ(\Omega,V)$ satisfies
\begin{equation}\label{prule}
dJ(\Omega,V)=I(\Omega)dH(\Omega,V)+H(\Omega)dI(\Omega,V).
\end{equation}
\end{lemma}

\begin{proof}
It follows from Definition \ref{defeulerder} of the Eulerian derivative.
\end{proof}
We want to point out, that we use indistinctly the notations $\frac{dJ(\Omega_{t})}{dt}|_{t=0}$ and $dJ(\Omega,V)$ to make reference to the Eulerian derivative at $\Omega$ in the direction of a vector field $V$. 

Now that we know that the product rule is holds in the sense of Eulerian differentiation, we proceed to find and expression for the Eulerian derivative of the shape functional given by the normalized first nontrivial Wentzell eigenvalue. 

Consider the eigenpair $(\lambda(t),u_{t})$ satisfying the boundary value problem
\begin{align}\label{wentzellperi}
     -\Delta u_{t} &=0,\quad x\in\Omega_{t}\\ \nonumber
     |\partial\Omega_{t}|(-\beta\Delta_{\tau}u_{t}+\partial_{n}u_{t}) &= \lambda(t)|\partial\Omega_{t}|u_{t}, \quad x\in\partial\Omega_{t}.\label{equ1}
\end{align} 
Multiplying the boundary condition by a test function $\phi_{t}$ such that $\partial_{n}\phi_{t}=0$, (see \cite{Dambrine}, \cite{Dambrine1}) and integrating by parts over the boundary $\partial\Omega_{t}$, we get the following weak formulation
\begin{equation}\label{weakformulation}
0=\int_{\partial\Omega_{t}}|\partial\Omega_{t}|\beta\nabla_{\tau}u_{t}\cdot\nabla_{\tau}\phi_{t}d\sigma_{t}+\int_{\partial\Omega_{t}}|\partial\Omega_{t}|\partial_{n}u_{t}\phi_{t}d\sigma_{t}-\int_{\partial\Omega_{t}}|\partial\Omega_{t}|\lambda(t)u_{t}\phi_{t}d\sigma_{t}.
\end{equation}
Taking the derivative with respect to $t$ in equation \eqref{weakformulation} and evaluating it at $t=0$ we have the following terms (see also, \cite{Caubet1})
\begin{equation}
\begin{aligned}
\RN{1}&=\frac{d}{dt}\left(\int_{\partial\Omega_{t}}|\partial\Omega_{t}|\beta\nabla_{\tau}u_{t}\cdot\nabla_{\tau}\phi_{t}d\sigma_{t}\right)_{t=0}\\
&= \frac{d(|\partial\Omega_{t}|)}{dt}\Big|_{t=0}\beta\int_{\partial\Omega}\nabla_{\tau}u\cdot\nabla_{\tau}\phi d\sigma+|\partial\Omega|\frac{d}{dt}\left(\int_{\partial\Omega_{t}}\beta\nabla_{\tau}u_{t}\cdot\nabla_{\tau}\phi_{t}d\sigma_{t}\right)_{t=0},\\
\RN{2}&=\frac{d}{dt}\left(\int_{\partial\Omega_{t}}|\partial\Omega_{t}|\partial_{n}u_{t}\phi_{t}d\sigma_{t}\right)_{t=0}\\
&= \frac{d(|\partial\Omega_{t}|)}{dt}\Big|_{t=0} \int_{\partial\Omega}\partial_{n}u\phi d\sigma 
+|\partial\Omega|\frac{d}{dt}\left(\int_{\partial\Omega_{t}}\partial_{n}u_{t}\phi_{t}\sigma_{t}\right)_{t=0},\\
\RN{3}&=\frac{d}{dt}\left(\int_{\partial\Omega_{t}}|\partial\Omega_{t}|\lambda(t)u_{t}\phi_{t}d\sigma_{t}\right)_{t=0}\\
&=(|\partial\Omega_{t}|\lambda(t))^{\prime}_{t=0}\int_{\partial\Omega}u\phi d\sigma + |\partial\Omega|\lambda \frac{d}{dt}\left(\int_{\partial\Omega_{t}}u_{t}\phi_{t}d\sigma_{t}\right)_{t=0}.
\end{aligned}
\end{equation}

From formula \eqref{shapederfor} we have that the term $$\frac{d}{dt}\left(\int_{\partial\Omega_{t}}\beta\nabla_{\tau}u_{t}\cdot\nabla_{\tau}\phi_{t}d\sigma_{t}\right)_{t=0}$$ in $\RN{1}$ becomes 
\begin{align}
\frac{d}{dt}\left(\int_{\partial\Omega_{t}}\beta\nabla_{\tau}u_{t}\cdot\nabla_{\tau}\phi_{t}d\sigma_{t}\right)_{t=0} &= \int_{\partial\Omega}\left[\beta V_{n}\frac{d(\nabla_{\tau}u_{t})}{dt}\Big|_{t=0}\cdot\nabla_{\tau}\phi+\beta\nabla_{\tau}u\cdot\frac{d(\nabla_{\tau}\phi)}{dt}\Big|_{t=0}\right]d\sigma \\ \nonumber
&+ \int_{\partial\Omega}V_{n}\left[\partial_{n}(\beta\nabla_{\tau}u\cdot\nabla_{\tau}\phi)+H(\beta\nabla_{\tau}u\cdot\nabla_{\tau}\phi)\right]d\sigma\\
&=\int_{\partial\Omega}\beta\left[\nabla_{\tau}u^{\prime}+\partial_{n}u\nabla_{\tau}V_{n}+(\nabla u\cdot\nabla_{\tau}V_{n})\vec{n}\right]\cdot\nabla_{\tau}\phi d\sigma\\ \nonumber
&+\int_{\partial\Omega}\beta\nabla_{\tau}u\cdot\left[\nabla_{\tau}\phi^{\prime}+(\nabla\phi\cdot\nabla_{\tau}V_{n})\vec{n}+\partial_{n}\phi\nabla_{\tau}V_{n}\right]d\sigma\\ \nonumber
&+\int_{\partial\Omega}V_{n}\left[\partial_{n}(\beta\nabla_{\tau}u\cdot\nabla_{\tau}\phi)+H(\beta\nabla_{\tau}u\cdot\nabla_{\tau}\phi)\right]d\sigma. \nonumber
\end{align}
Where we have used the fact that $\frac{d(\nabla_{\tau}u)}{dt}\Big|_{t=0}=\nabla_{\tau}u^{\prime}+(\nabla u\cdot \nabla_{\tau}V_{n})\vec{n}+\partial_{n}u\nabla_{\tau}V_{n}$ (see proof of Proposition 5.1 in \cite{Caubet1}) to go from the first equality to the second. Now, $\nabla_{\tau}\phi\cdot\vec
n=\nabla_{\tau}u\cdot\vec{n}=0$, $\phi^{\prime}(\partial\Omega;V)=\phi^{\prime}(\partial\Omega;V(0))=\dot{\phi}(\partial\Omega;V(0))-\nabla_{\tau}\phi(\partial\Omega)\cdot V(0)=0$ (Here we are using the fact that $\dot{\phi}(\partial\Omega;V(0))=\nabla_{\tau}\phi(\partial\Omega)\cdot V(0)$, see Proposition 2.77 page 101 and Proposition 2.89 page 114 in \cite{sokolowski}) and we also have $\phi^{\prime}(\partial\Omega;V)=\dot{\phi}(\partial\Omega;V(0))-\left\langle \nabla\phi(\partial\Omega),V(0)\right\rangle_{\mathbb{R}^{N}}+\partial_{n}\phi\left\langle V(0),\vec{n}\right\rangle_{\mathbb{R}^{N}}$ which implies $\phi^{\prime}(\partial\Omega;V)=\partial_{n}\phi\left\langle V(0),\vec{n}\right\rangle_{\mathbb{R}^{N}}=0$ and hence $\nabla_{\tau}\phi^{\prime}=0$. The dot notation $\dot{\phi}(\partial\Omega,V)$ is the material derivative of $\phi(\partial\Omega)$, in the direction of the vector field $V$ \cite{sokolowski}. Thus, 
\begin{align}
&\frac{d}{dt}\left(\int_{\partial\Omega_{t}}\beta\nabla_{\tau}u_{t}\cdot\nabla_{\tau}\phi_{t}d\sigma_{t}\right)_{t=0} = \int_{\partial\Omega}\beta\left[\nabla_{\tau}u^{\prime}+\partial_{n}u\nabla_{\tau}V_{n}\right]\cdot\nabla_{\tau}\phi d\sigma\\ \nonumber
&+ \int_{\partial\Omega}V_{n}\left[\partial_{n}(\beta\nabla_{\tau}u\cdot\nabla_{\tau}\phi)+H(\beta\nabla_{\tau}u\cdot\nabla_{\tau}\phi)\right]d\sigma \\ \label{510}
&= \int_{\partial\Omega}\beta\nabla_{\tau}u^{\prime}\cdot\nabla_{\tau}\phi d\sigma+\int_{\partial\Omega}\beta\partial_{n}u\nabla_{\tau}V_{n}\cdot\nabla_{\tau}\phi d\sigma \\ \nonumber
&+\int_{\partial\Omega}V_{n}\left[\partial_{n}(\beta\nabla_{\tau}u\cdot\nabla_{\tau}\phi)+H(\beta\nabla_{\tau}u\cdot\nabla_{\tau}\phi)\right]d\sigma.  \label{equations}
\end{align}

The first integral term in \eqref{510} can be written as 
$$\int_{\partial\Omega}\beta\nabla_{\tau}u^{\prime}\cdot\nabla_{\tau}\phi d\sigma=\int_{\partial\Omega}-\beta\Delta_{\tau}u^{\prime}\phi d\sigma.$$
Since $\partial_{n}\phi=0$, using Theorem B.2 in \cite{Caubet1} we obtain
\begin{align}
\int_{\partial\Omega}V_{n}[\partial_{n}(\beta\nabla_{\tau}u\cdot\nabla_{\tau}\phi)d\sigma]=\beta\int_{\partial\Omega}\left(\nabla_{\tau}(\partial_{n}u)\cdot\nabla_{\tau}\phi+\nabla_{\tau}(\partial_{n}\phi)\cdot\nabla_{\tau}u-2D^{2}b\nabla_{\tau}u\cdot\nabla_{\tau}\phi\right)V_{n}d\sigma,
\end{align}
so,
\begin{align}
\frac{d}{dt}\left(\int_{\partial\Omega_{t}}\beta\nabla_{\tau}u_{t}\cdot\nabla_{\tau}\phi_{t}d\sigma_{t}\right)_{t=0} = 
\int_{\partial\Omega}[-\beta\Delta_{\tau}u^{\prime}-\beta\text{div}_{\tau}(V_{n}(HId-2D^{2}b)\nabla_{\tau}u)-\beta\Delta_{\tau}(V_{n}\partial_{n}u)]\phi d\sigma.
\end{align}

Now we focus in the second term of $\RN{2}$, using \eqref{shapederfor} again we have 
\begin{align}
\frac{d}{dt}\left(\int_{\partial\Omega_{t}}\partial_{n}u_{t}\phi_{t}d\sigma_{t}\right)_{t=0} &= 
\integral (\partial_{n}u\phi)^{\prime} d\sigma + \integral V_{n}(\partial_{n}(\partial_{n}u\phi)+H\partial_{n}u\phi) d\sigma. \label{523}
\end{align}
Based on our choice of $\phi$, we have $(\partial_{n}u\phi)^{\prime}_{t=0}=(\partial_{n}u)^{\prime}\phi$, since $(\partial_{n}u)^{\prime}=\nabla u^{\prime}\cdot n+\nabla u \cdot (\vec{n})^{\prime}$.
Using that $(\vec{n})^{\prime}=-\nabla_{\tau}V_{n}$ (see Prop. B.1, \cite{Caubet1} ). Equation \eqref{523} becomes
\begin{align}\label{524}
\frac{d}{dt}\left(\int_{\partial\Omega_{t}}\partial_{n}u_{t}\phi_{t}d\sigma_{t}\right)_{t=0} &= 
\integral \partial_{n}u^{\prime}\phi d\sigma-\integral \nabla_{\tau}V_{n}\cdot\nablat u \phi d\sigma + \integral V_{n}[(\Delta u - \Delta_{\tau}u)\phi+\partial_{n}u\partial_{n}\phi]d\sigma.
\end{align}
Equation \eqref{524} was obtained using the relation $\Delta u - \Delta_{\tau}u = H\partial_{n}u+\partial^{2}_{nn}u$ and the identity $\nablat u\cdot \nablat V_{n}=\nabla u \cdot \nablat V_{n}$. 
Similarly, using formula \eqref{shapederfor} we get that the second term in $\RN{3}$ can be written as
\begin{align}
|\partial\Omega|\lambda(0)\integral u^{\prime}\phi d\sigma +|\partial\Omega|\lambda(0)\integral V_{n}(\partial_{n}u\phi+\partial_{n}\phi+Hu\phi)d\sigma.
\end{align}

Thus, $\RN{1}$, $\RN{2}$ and $\RN{3}$ can be written respectively as:
\begin{align}
\nonumber
&\frac{d}{dt}\left(\int_{\partial\Omega_{t}}|\partial\Omega_{t}|\beta\nabla_{\tau}u_{t}\cdot\nabla_{\tau}\phi_{t}d\sigma_{t}\right)_{t=0} = \frac{d(|\partial\Omega_{t}|)}{dt}\Big|_{t=0}\beta\int_{\partial\Omega}\nabla_{\tau}u\cdot\nabla_{\tau}\phi d\sigma \\ \label{526}
&+|\partial\Omega|\integral [-\beta\Delta_{\tau}u^{\prime}-\beta \text{div}_{\tau}(V_{n}(HId-2D^{2}b)\nabla_{\tau}u)-\beta \Delta_{\tau}(V_{n}\partial_{n}u)]\phi d\sigma 
\end{align}

\begin{align}
\nonumber
\frac{d}{dt}\left(\int_{\partial\Omega_{t}}|\partial\Omega_{t}|\partial_{n}u_{t}\phi_{t}d\sigma_{t}\right)_{t=0}&= \frac{d(|\partial\Omega_{t}|)}{dt}\Big|_{t=0} \int_{\partial\Omega}\partial_{n}u\phi d\sigma \\ \label{527}
&+|\partial\Omega|\left(\integral(\partial_{n}u^{\prime}\phi d\sigma - \text{div}_{\tau}(V_{n}\nablat u)\phi d\sigma\right)
\end{align}

\begin{align}
\nonumber
\frac{d}{dt}\left(\int_{\partial\Omega_{t}}|\partial\Omega_{t}|\lambda(t)u_{t}\phi_{t}d\sigma_{t}\right)_{t=0}&=(|\partial\Omega_{t}|\lambda(t))^{\prime}_{t=0}\int_{\partial\Omega}u\phi d\sigma \\
&+ |\partial\Omega|\lambda(0)\integral u^{\prime}\phi d\sigma +  |\partial\Omega|\lambda(0)\integral V_{n}(\partial_{n}u+Hu)\phi d\sigma. \label{528}
\end{align}
Where in \eqref{527} we have used the fact that $\Delta u = 0$ and $\text{div}_{\tau}(V_{n}\nablat u)=\nablat V_{n}\cdot\nablat u + V_{n}\text{div}_{\tau}(\nablat u)$ (see, Pag 91 in \cite{sokolowski}). Then, from \eqref{526}, \eqref{527} and \eqref{528} we get the following expression
\begin{align}
\nonumber
\frac{d(|\partial\Omega_{t}|)}{dt}\Big|_{t=0}\beta\int_{\partial\Omega}\nabla_{\tau}u\cdot\nabla_{\tau}\phi d\sigma &+|\partial\Omega|\integral [-\beta\Delta_{\tau}u^{\prime}-\beta \text{div}_{\tau}(V_{n}(HId-2D^{2}b)\nabla_{\tau}u)-\beta \Delta_{\tau}(V_{n}\partial_{n}u)]\phi d\sigma \\ \nonumber
&+\frac{d(|\partial\Omega_{t}|)}{dt}\Big|_{t=0} \int_{\partial\Omega}\partial_{n}u\phi d\sigma +|\partial\Omega|\integral(\partial_{n}u^{\prime}- \text{div}_{\tau}(V_{n}\nablat u))\phi d\sigma \\ \label{529}
&= (|\partial\Omega_{t}|\lambda(t))^{\prime}_{t=0}\int_{\partial\Omega}u\phi d\sigma + |\partial\Omega|\lambda(0)\integral u^{\prime}\phi d\sigma\\ \nonumber 
&+  |\partial\Omega|\lambda(0)\integral V_{n}(\partial_{n}u+Hu)\phi d\sigma. 
\end{align}

Reordering \eqref{529} we have
\begin{align}
\nonumber
&|\partial\Omega|\integral [-\beta\Delta_{\tau}u^{\prime}-\beta \text{div}_{\tau}(V_{n}(HId-2D^{2}b)\nabla_{\tau}u)-\beta \Delta_{\tau}(V_{n}\partial_{n}u)]\phi d\sigma \\ \label{530}
&+ |\partial\Omega|\integral[\partial_{n}u^{\prime}-\text{div}_{\tau}(V_{n}\nablat u)]\phi d\sigma - (|\partial\Omega|\lambda)\integral u^{\prime}\phi d\sigma \\ \nonumber
&= \frac{d(|\partial\Omega_{t}|)}{dt}\Big|_{t=0}\integral
\beta\Delta_{\tau}u\phi d\sigma - \frac{d(|\partial\Omega_{t}|)}{dt}\Big|_{t=0}\integral\partial_{n}u\phi d\sigma + \integral(|\partial\Omega_{t}|\lambda(t))^{\prime}_{t=0}u\phi d\sigma \\ \nonumber
&+|\partial\Omega|\lambda\integral V_{n}(\partial_{n}u+Hu)\phi d\sigma.
\end{align}
Notice that equation \eqref{530} implies the following expression, analogous to the one given by Theorem \ref{thm2}. 
\begin{align}
\nonumber
|\partial\Omega|(-\beta\Delta_{\tau}u^{\prime}+\partial_{n}u^{\prime}-\lambda u^{\prime}) &= \beta|\partial\Omega|\Delta_{\tau}(V_{n}\partial_{n}u)-\beta|\partial\Omega|\text{div}_{\tau}(V_{n}(2D^{2}b-HId)\nablat u) \\ \nonumber
&+\text{div}_{\tau}(V_{n}\nablat u)|\partial\Omega|+(|\partial\Omega_{t}|\lambda(t))^{\prime}_{t=0}u+\lambda|\partial\Omega|V_{n}(\partial_{n}u+Hu)\\ 
&+K(V)(\beta\Delta_{\tau}u-\partial_{n}u).
\end{align}
Where $K(V)$ is the Eulerian derivative of the perimeter at time $t=0$ i.e. (see, page 55 in \cite{sokolowski})
$$K(V)=\frac{d}{dt}\left(\integral d\sigma_{t}\right)_{t=0}=\integral \text{div}_{\tau}(V(0))d\sigma.$$
Then, using the Wentzell boundary condition we obtain.
\begin{align}
\nonumber
|\partial\Omega|(-\beta\Delta_{\tau}u^{\prime}+\partial_{n}u^{\prime}-\lambda u^{\prime}) &= \beta|\partial\Omega|\Delta_{\tau}(V_{n}\partial_{n}u)-\beta|\partial\Omega|\text{div}_{\tau}(V_{n}(2D^{2}b-HId)\nablat u) \\ \nonumber
&+\text{div}_{\tau}(V_{n}\nablat u)|\partial\Omega|+(|\partial\Omega_{t}|\lambda(t))^{\prime}_{t=0}u+\lambda|\partial\Omega|V_{n}(\partial_{n}u+Hu)\\ \label{532}
&-K(V)\lambda u.
\end{align}

Following the proof for the derivative of a simple eigenvalue in \cite{Dambrine}, if we multiply both sides of the equality in \eqref{532} by the normalized eigenfunction $u$ and integrate over $\partial\Omega$, we obtain that the left-hand side gives:
\begin{align}
|\partial\Omega|\integral(-\beta\Delta_{\tau}u^{\prime}+\partial_{n}u^{\prime}-\lambda u^{\prime})u d\sigma &= |\partial\Omega|\integral(-\beta\Delta_{\tau}u+\partial_{n}u-\lambda u)u^{\prime}d\sigma = 0. \label{533}
\end{align}
where the first equality follows from Green's identities and the last equality is just the boundary condition. On the other hand, the terms in the right hand side give
\begin{align}
|\partial\Omega|\integral \beta u \Delta_{\tau}(V_{n}\partial_{n}u)d\sigma &= |\partial\Omega|\integral\beta\Delta_{\tau}u(V_{n}\partial_{n}u)d\sigma = |\partial\Omega|\integral(\partial_{n}u-\lambda u)(V_{n}\partial_{n}u)d\sigma, \label{534}
\end{align}

\begin{align}
|\partial\Omega|\integral\beta \text{div}_{\tau}(V_{n}(2D^{2}b-HId)\nablat u)u d\sigma &= -|\partial\Omega|\integral
\beta V_{n}(2D^{2}b-HId)\nablat u \cdot \nablat u d\sigma, \label{535}
\end{align}

\begin{align}
|\partial\Omega|\integral \text{div}_{\tau}(V_{n}\nablat u)u d\sigma = -|\partial\Omega|\integral V_{n}|\nablat u|^{2}d\sigma, \label{536}
\end{align}

\begin{align}
\integral (\lambda(t)|\partial\Omega|)^{\prime}u^{2} d\sigma = (\lambda(t)|\partial\Omega_{t}|)^{\prime}_{t=0}, \label{537}
\end{align}

\begin{align}
\integral K(V)(-\lambda u^{2})d\sigma =-\lambda(0)K(V). \label{538}
\end{align}
In equations \eqref{537}, \eqref{538} we have used the normalization condition $\int_{\partial\Omega}u^{2}d\sigma=1$ also used in the proof of Theorem \ref{simplecase} (see \cite{Dambrine}).
Then from \eqref{534}, \eqref{535}, \eqref{536}, \eqref{537}, \eqref{538} we have
\begin{align}
\nonumber
0&=|\partial\Omega|\Bigg(\integral [V_{n}(\partial_{n}u)^{2}-\lambda u V_{n}\partial_{n}u+V_{n}\beta(2D^2{b}-HId)\nablat u \cdot \nablat u-V_{n}|\nablat u|^{2} \\
&+\lambda u V_{n}\partial_{n}u+\lambda V_{n}H|u|^{2}]d\sigma\Bigg) + (\lambda(t)|\partial\Omega_{t}|)^{\prime}_{t=0} - \lambda(0)K(V). \label{539}
\end{align}
Therefore, we have proved the following proposition.

\begin{prop}\label{propo1}
If $\lambda$ is a simple eigenvalue of the Wentzell problem, the map $t\mapsto |\partial\Omega_{t}|\lambda(t)$ has derivative at $t=0$ given by 
\begin{align}
\nonumber
(\lambda(t)|\partial\Omega_{t}|)^{\prime}_{t=0} &= \lambda(0)K(V) + |\partial\Omega|\Bigg(\integral V_{n}[-(\partial_{n}u)^{2}+\beta(HId-2D^2{b})\nablat u \cdot \nablat u \\
&+|\nablat u|^{2}-\lambda H|u|^{2}]d\sigma\Bigg) \\ \nonumber
&=\left(\lambda(t)\frac{d|\partial\Omega_{t}|}{dt}+|\partial\Omega_{t}|\lambda^{\prime}(t)\right)_{t=0}
\end{align}
where $\lambda^{\prime}(0)$ is given by Theorem \ref{simplecase} in Section \ref{secshapecal}.
\end{prop}
Since in our case we are dealing with a multiple eigenvalue, Proposition \ref{propo1} is not the one that we are interested in. However, we can modify it to get an analogous result for the case of a multiple eigenvalue. As in the simple case, the computations and ideas follow those found in \cite{Dambrine1} and \cite{Caubet1}. 

We begin by considering a smooth branch $t\mapsto (u(t,x),\lambda(\Omega_{t}))$ corresponding to a multiple eigenvalue $\lambda$ of the Wentzell eigenvalue problem (this branch in itself is simple so it satisfies the boundary condition of the boundary value problems in Theorem \ref{thm2} in Section \ref{secshapecal}. In our case, it satisfies the boundary condition
\begin{align}
\nonumber
|\partial\Omega|(-\beta\Delta_{\tau}u^{\prime}+\partial_{n}u^{\prime}-\lambda u^{\prime})&= |\partial\Omega|(\beta\Delta_{\tau}(V_{n}\partial_{n}u)-\beta \text{div}_{\tau}(V_{n}(2D^{2}b-HId)\nablat u)+\text{div}_{\tau}(V_{n}\nablat u)\\ \label{541}
&+\lambda V_{n}(\partial_{n}u+Hu))+(\lambda(t)|\partial\Omega_{t}|)^{\prime}_{t=0}u+K(V)(\beta\Delta_{\tau}u-\partial_{n}u).
\end{align}

Using the decomposition $u=\sum_{i=1}^{m}d_{i}u_{i}$, where $u_{i}$ are the normalized eigenfunctions of the eigenspace associated to $\lambda$, multiplying both sides of \eqref{541} by $u_{j}$ and integrating over $\partial\Omega$ we get that the left hand side of \eqref{541} satisfies
\begin{align}
|\partial\Omega|\left(-\beta\Delta_{\tau}\left(\sum_{i}d_{i}u_{i}^{\prime}\right)+\partial_{n}\left(\sum_{i}d_{i}u_{i}^{\prime}\right)-\lambda\sum_{i}d_{i}u_{i}^{\prime}\right) = |\partial\Omega|\sum_{i}d_{i}\left(-\beta\Delta_{\tau}u_{i}^{\prime}+\partial_{n}u_{i}^{\prime}-\lambda u_{i}^{\prime}\right).
\end{align}
Therefore
\begin{align}
|\partial\Omega|\integral u_{j}\left(-\beta\Delta_{\tau}u_{i}^{\prime}+\partial_{n}u_{i}^{\prime}-\lambda u_{i}^{\prime}\right)d\sigma = |\partial\Omega|\integral u_{j}^{\prime}\left(-\beta\Delta_{\tau}u_{i}+\partial_{n}u_{i}-\lambda u_{i}\right)d\sigma = 0.\label{543}
\end{align}

On the other hand, working term by term with the right hand side of \eqref{541} side we get
\begin{align}
\nonumber
\integral \beta\Delta_{\tau}\left(V_{n}\partial_{n}\left(\sum_{i}d_{i}u_{i}\right)\right)u_{j}d\sigma &= \beta\sum_{i}d_{i}\integral\Delta_{\tau}(V_{n}\partial_{n}u_{i})u_{j}d\sigma \\ \label{544}
&=-\beta\sum_{i}d_{i}\integral\nablat u_{j}\cdot\nablat(V_{n}\partial_{n}u_{i})d\sigma \\ \nonumber
&=\beta\sum_{i}d_{i}\integral\Delta_{\tau}u_{j}(V_{n}\partial_{n}u_{i})d\sigma \\ \nonumber
&=\sum_{i}d_{i}\integral (\partial_{n}u_{j}-\lambda u_{j})(V_{n}\partial_{n}u_{i})d\sigma, 
\end{align}

\begin{align}
\nonumber
 -\beta\integral \text{div}_{\tau}\left(V_{n}(2D^{2}b-HId)\nablat\left(\sum_{i}d_{i}u_{i}\right)\right)u_{j}d\sigma &= \integral \beta\nablat u_{j}\cdot [(V_{n}(2D^{2}b-HId))\sum_{i}d_{i}\nablat u_{i}]d\sigma \\ \label{545}
 &= \sum_{i}d_{i}\integral \beta (V_{n}(2D^{2}b-HId))\nablat u_{i}\cdot \nablat u_{j}d\sigma,
\end{align}

\begin{align}
\integral \text{div}_{\tau}\left(V_{n}\nablat \left(\sum_{i}d_{i}u_{i}\right)\right)u_{j}d\sigma &= \sum_{i}d_{i}\left(-\integral V_{n}\nablat u_{j}\cdot \nablat u_{i}d\sigma\right), \label{546}
\end{align}

\begin{align}
\integral \lambda V_{n}\left(\partial_{n}\left(\sum_{i}d_{i}u_{i}\right) +H\left(\sum_{i}d_{i} u_{i}\right)\right)u_{j}d\sigma &= \sum_{i}d_{i}\integral V_{n}(\lambda u_{j}\partial_{n}u_{i}+\lambda Hu_{j}u_{i})d\sigma, \label{547}
\end{align}

\begin{align}
\integral (|\partial\Omega_{t}|\lambda(t))^{\prime}_{t=0}\left(\sum_{i}d_{i}u_{i}\right)u_{j}d\sigma = (|\partial\Omega_{t}|\lambda(t))^{\prime}_{t=0}\sum_{i}d_{j}\integral u_{i}u_{j}d\sigma, \label{548}
\end{align}

\begin{align}
\integral K(V)\left(-\lambda\sum_{i}d_{i}u_{i}u_{j}\right)d\sigma = -K(V)\lambda\sum_{i}d_{i}\integral u_{i}u_{j}d\sigma. \label{549}
\end{align}

Then, using the expression for $u^{\prime}$ obtained in \eqref{541} and the expressions \eqref{543},  \eqref{544}, \eqref{545}, \eqref{546}, \eqref{547}, \eqref{548} and \eqref{549} we obtain 
\begin{align}
&|\partial\Omega|\Bigg[\sum_{i}d_{i}\Bigg(\integral [(\partial_{n}u_{j}-\lambda u_{j})V_{n}\partial_{n}u_{i}+\beta(V_{n}(2D^{2}b-HId))\nablat u_{j}\cdot \nablat u_{i}]d\sigma \\ \nonumber
&-\integral V_{n}\nablat u_{i}\cdot \nablat u_{j}d\sigma+\integral V_{n}(\lambda u_{j}\partial_{n}u_{i}+\lambda H u_{i}u_{j})d\sigma\Bigg)\Bigg]+\sum_{i}d_{i}(|\partial\Omega_{t}|\lambda(t))^{\prime}_{t=0}\integral u_{i}u_{j}d\sigma\\ \nonumber
&- \sum_{i}d_{i}K(V)\lambda\integral u_{i}u_{j}d\sigma = 0,
\end{align}
which we can rewrite as
\begin{align}
\nonumber
&|\partial\Omega|\Bigg[\sum_{i}d_{i}\integral V_{n}(\partial_{n}u_{j}\partial_{n}u_{i}+\beta(V_{n}(2D^{2}b-HId))\nablat u_{j}\cdot \nablat u_{i}-\nablat u_{j}\cdot \nablat u_{i}+\lambda Hu_{i}u_{j})d\sigma\Bigg]\\ \label{551}
&+\sum_{i}d_{i}\integral u_{i}u_{j}d\sigma\left((|\partial\Omega_{t}|\lambda(t\right))^{\prime}_{t=0}-K(V)\lambda)=0.
\end{align}

Setting $N_{ij}=\integral u_{i}u_{j}d\sigma=\delta_{ij}$ we can write \eqref{551} in the form
\begin{align}
&|\partial\Omega|\sum_{i}d_{i}(-M_{ij})+\sum_{i}d_{i}\delta_{ij}[(|\partial\Omega_{t}|\lambda(t))^{\prime}_{t=0}-K(V)\lambda] = 0, \label{552}
\end{align}
and therefore
\begin{align}
\sum_{i}d_{i}\delta_{ij}(|\partial\Omega_{t}|\lambda(t))^{\prime}_{t=0} &= |\partial\Omega|\sum_{i}d_{i}M_{ij}+K(V)\lambda\sum_{i}d_{i}\delta_{ij},
\end{align}
where $\delta_{ij}$ is the Kronecker delta. Noting $M_{ij}$ is the matrix defined in Theorem \ref{multievals}, section \ref{secshapecal}. We have proved the following result (compare with Thm. \ref{multievals}, in Sec. \ref{secshapecal}).

\begin{prop}

Let $\lambda$ be a multiple eigenvalue of order $m \geq 2$. Then each $t \mapsto \lambda_{i}(t)|\partial\Omega_{t}|$ for $i=1\dots N$ has a derivative near $0$, and the values of $(\lambda_{i}(t)|\partial\Omega_{t}|)^{\prime}_{t=0}$, $i=1\dots N$ satisfy 
\begin{align}
(\lambda(t)|\partial\Omega_{t}|)^{\prime}_{t=0}\vec{d}=(|\partial\Omega|M+K(V)\lambda Id)\vec{d},\label{554}
\end{align}
where $M$ is the matrix defined by Theorem \ref{multievals}, in Section \ref{secshapecal}, $K(V)=\frac{d|\partial\Omega_{t}|}{dt}\Big|_{t=0}$, and $\vec{d}$ is the vector of coefficients in the decomposition of $u(0,x)=\sum_{i=1}^{m}d_{i}u_{i}$, where $u_{i}$ are the corresponding basis for the eigenspace. Here $^{\prime}$ denotes the shape derivative in the direction of $V$.
\end{prop}

That is, the derivatives of $(\lambda_{i}(t)|\partial\Omega_{t}|)^{\prime}_{t=0}$ are the eigenvalues of the matrix $$|\partial\Omega|M+K(V)\lambda Id.$$ 
In other words, the product rule holds in the sense described by equation \eqref{554}. That is 
\begin{align*}
(\lambda(t)|\partial\Omega_{t}|)^{\prime}_{t=0}\vec{d}&=|\partial\Omega|M\vec{d}+K(V)\lambda\vec{d} \\
&= |\partial\Omega|\lambda^{\prime}(0)\vec{d}+K(V)\lambda\vec{d}.
\end{align*}

\subsection{The annulus case}
For the case of an annulus, solutions to $\Delta u=0$ are given by $f_{k}(r,\theta)=(A_{k}r^{k}+A_{-k}r^{-k})\cos(k\theta)$ or  $f_{k}(r,\theta)=(A_{k}r^{k}+A_{-k}r^{-k})\sin(k\theta)$. Since we want $f_{k}$ to be the eigenfunction associated to the corresponding eingenvalue $\lambda_{k}$, the constants $A_{k}, A_{-k}$ satisfy the following system:
\begin{align*}
A_{k}(\beta k^{2}+k-\lambda_{k})+A_{-k}(\beta k^{2}-k-\lambda_{k})&=0\\
A_{k}(\beta k^{2}\epsilon^{k-2}-k\epsilon^{k-1}-\lambda_{k}\epsilon^{k})+A_{-k}(\beta k^{2}\epsilon^{-k-2}+k\epsilon^{-k-1}-\lambda_{k}\epsilon^{-k})&=0.
\end{align*}

We consider an annulus $\Omega_{\epsilon}=\mathbb{D}\setminus B(0,\epsilon)$, with $\epsilon \in (0,1)$. For this domain $\Omega_{\epsilon}$, the signed distance $b$ is defined as
\begin{align}
b(r,\theta)=
\begin{cases} 
      r-1, & r\geq \frac{1+\epsilon}{2}, \\
      \epsilon-r,  & 0 \leq r < \frac{1+\epsilon}{2}.
\end{cases}
\end{align}
From this it follows that 
\begin{align}
\nabla b(r,\theta)|_{\partial\Omega_{\epsilon}}=
\begin{cases} 
      \left\langle \cos\theta,\sin\theta\right\rangle, & r=1,\\
      \left\langle -\cos\theta,-\sin\theta\right\rangle, & r=\epsilon.
\end{cases}
\end{align}
Then, we have that on $\partial\mathbb{D}$, $D^{2}b(r,\theta)=\frac{Id}{r}-B$ where $B$ is given by
\begin{align}
B=
\begin{bmatrix}
    \frac{\cos^{2}\theta }{r}     & \frac{r\sin\theta r\cos\theta}{r^3} \\
    \frac{r\sin\theta r\cos\theta}{r^3}       & \frac{\sin^{2}\theta }{r} 
\end{bmatrix}.
\end{align}
On the other hand, we have $\nablat u_{i}=\frac{1}{r}\frac{\partial u_{i}}{\partial\theta}\vec{n}_{\theta}$, and $\nablat u_{i}\cdot\nablat u_{j}=\frac{1}{r^2}\frac{\partial u_{i}}{\partial\theta}\frac{\partial u_{j}}{\partial\theta}$. Thus, 
$$(HId-2D^{2}b)\nablat u_{i}\cdot\nablat u_{j}=-\frac{1}{r^3}\frac{\partial u_{i}}{\partial\theta}\frac{\partial u_{j}}{\partial\theta}.$$

Since we are interested in the Steklov eigenvalues, setting $\beta=0$ we get that the entries of the matrix defined in Theorem \ref{multievals} on $\partial\mathbb{D}$ are given by
\begin{align}
\widetilde{M}_{ij}=\int_{0}^{2\pi}V_{n}\left(\frac{-(A_{1}+A_{-1})^{2}-(A_{1}-A_{1})^{2}-\lambda_{1}(A_{1}+A_{-1})^{2}}{\|g_{1}\|\|g_{2}\|}\right)\cos\theta\sin\theta d\theta, \quad \text{if} \quad i\neq j
\end{align}
\begin{align}
\widetilde{M}_{11}=\int_{0}^{2\pi}V_{n}\left(\frac{(A_{1}+A_{-1})^{2}}{\|g_{1}\|^{2}}(\sin^{2}\theta-\lambda_{1}\cos^{2}\theta)-\frac{(A_{1}-A_{-1})^{2}}{\|g_{1}\|^{2}}\cos^{2}\theta\right) d\theta,
\end{align}
\begin{align}
\widetilde{M}_{22}=\int_{0}^{2\pi}V_{n}\left(\frac{(A_{1}+A_{-1})^{2}}{\|g_{2}\|^{2}}(\cos^{2}\theta-\lambda_{1}\sin^{2}\theta)-\frac{(A_{1}-A_{-1})^{2}}{\|g_{2}\|^{2}}\sin^{2}\theta\right) d\theta.
\end{align}
Similarly on $\partial B(0,\epsilon)$ we get
\begin{align}
M_{ij}=C(\epsilon,\lambda_{1})\int_{0}^{2\pi}V_{n}\cos\theta\sin\theta d\theta, \quad \text{if} \quad  i\neq j
\end{align}
where
\begin{align}
C(\epsilon,\lambda_{1})=\left(\frac{(A_{1}\epsilon+A_{-1}\epsilon^{-1})^{2}}{\|g_{1}\|\|g_{2}\|}\left(\frac{\lambda_{1}}{\epsilon}-\frac{1}{\epsilon^{2}}\right)-\frac{(A_{1}-A_{-1}\epsilon^{-2})^{2}}{\|g_{1}\|\|g_{2}\|}\right)\epsilon
\end{align}
and
\begin{align}
M_{11}=\int_{0}^{2\pi}V_{n}\left(\frac{(A_{1}\epsilon+A_{-1}\epsilon^{-1})^{2}}{\|g_{1}\|^{2}}\left(\frac{\sin^{2}\theta}{\epsilon^{2}}+\lambda_{1}\frac{\cos^{2}\theta}{\epsilon}\right)-\frac{(A_{1}-A_{-1}\epsilon^{-2})^{2}}{\|g_{1}\|^{2}}\cos^{2}\theta\right)\epsilon d\theta
\end{align}
\begin{align}
M_{22}=\int_{0}^{2\pi}V_{n}\left(\frac{(A_{1}\epsilon+A_{-1}\epsilon^{-1})^{2}}{\|g_{2}\|^{2}}\left(\frac{\cos^{2}\theta}{\epsilon^{2}}+\lambda_{1}\frac{\sin^{2}\theta}{\epsilon}\right)-\frac{(A_{1}-A_{-1}\epsilon^{-2})^{2}}{\|g_{2}\|^{2}}\sin^{2}\theta\right)\epsilon d\theta,
\end{align}
where $g_{1}(r,\theta)=(A_{1}r+A_{-1}r^{-1})\cos{\theta}$, $g_{2}(r,\theta)=(A_{1}r+A_{-1}r^{-1})\sin{\theta}$. Thus, the matrix defined in Theorem \ref{multievals} for the domain $\Omega_{\epsilon}$ is given by $S=-M+\widetilde{M}$, with entries $S_{ij}=-M_{ij}+\widetilde{M}_{ij}$. 

Now that we know the related quantities for the particular case of an annulus $\Omega_{\epsilon}$, we will consider perturbations of such domains. Namely, we fix the outer boundary while we ``deform'' the inner boundary. Thus, for this type of deformations, we consider a vector field $V \in  C((0,\epsilon);V^{k}(D))$ satisfying $V|_{\partial\mathbb{D}}=0$.
Put $\|g_{1}\|=\|g_{2}\|=G$ and set 
\begin{align}
C_{1}(\epsilon,\lambda_{1})&=\left(\frac{(A_{1}\epsilon+A_{-1}\epsilon^{-1})^{2}}{G^{2}}\left(\frac{\lambda_{1}}{\epsilon}-\frac{1}{\epsilon^{2}}\right)-\frac{(A_{1}-A_{-1}\epsilon^{-2})^{2}}{G^{2}}\right)\epsilon, \\
C_{2}(\epsilon)&=\frac{(A_{1}\epsilon+A_{-1}\epsilon^{-1})^{2}}{\epsilon G^{2}}=C_{5}(\epsilon),\\ 
C_{3}(\epsilon,\lambda_{1})&=\left(\frac{(A_{1}\epsilon+A_{-1}\epsilon^{-1})^{2}}{\epsilon G^{2}}\lambda_{1}-\frac{(A_{1}-A_{-1}\epsilon^{-2})^{2}}{G^{2}}\right)\epsilon=C_{4}(\epsilon,\lambda_{1}).
\end{align}
Since the condition  $V|_{\partial\mathbb{D}}=0$ makes $\widetilde{M}=0$, the  entries of the matrix $M$ can be writen in a more compact way as follows 
\begin{align}
M_{ij}&=C_{1}(\epsilon,\lambda_{1})\int_{0}^{2\pi}V_{n}\sin\theta\cos\theta d\theta,  \quad i\neq j, \\
M_{11}&=C_{2}(\epsilon)\int_{0}^{2\pi}V_{n}\sin^{2}\theta d\theta + C_{3}(\epsilon,\lambda_{1})\int_{0}^{2\pi}V_{n}\cos^{2}\theta d\theta, \\
M_{22}&=C_{4}(\epsilon,\lambda_{1})\int_{0}^{2\pi}V_{n}\sin^{2}\theta d\theta + C_{5}(\epsilon)\int_{0}^{2\pi}V_{n}\cos^{2}\theta d\theta.
\end{align}
Note that $C_{1}=C_{3}-C_{2}$.

In addition to the condition  $V|_{\partial\mathbb{D}}=0$, we also consider the following decomposition of the normal component $V_{n}$ of $V$ on $\partial B(0,\epsilon)$. We write $V_{n}(\theta)=\omega_{n}^{r}+\omega_{n}^{l}(\theta)$, where $\omega_{n}^{r}=k$ is the component that produces radial symmetric deformations of the inner boundary, and $\omega_{n}^{l}$ produces deformations that are length preserving. Note that this decomposition is possible since we can write $V_{n}$ as follows.
\begin{align}
V_{n}=\frac{1}{|\partial\Omega|}\int_{\partial\Omega} V_{n}+\left(V_{n}-\frac{1}{|\partial\Omega|}\int_{\partial\Omega} V_{n}\right).
\end{align}

Setting $\omega_{n}^{r}=\frac{1}{|\partial\Omega|}\int_{\partial\Omega} V_{n}$ and $\omega_{n}^{l}=V_{n}-\frac{1}{|\partial\Omega|}\int_{\partial\Omega} V_{n}$, we see that $\frac{1}{|\partial\Omega|}\int_{\partial\Omega} \omega_{n}^{l}=0$, thus $\omega_{n}^{l}$ indeed preserves lenghts. Notice that by writing $V_{n}$ in this form, we are considering deformations that preserve the length of the inner boundary, once a radial deformation has been applied at a given ``instant''.

Given that $\lambda_{1}$ has multiplicity two, we have from Theorem \ref{thm1}, that there exists two maps
\begin{align}
t\mapsto \lambda_{1,1}(t), \quad t\mapsto\lambda_{1,2}(t),
\end{align}
that are differentiable in a neighborhood $0$. Using the decomposition  $V_{n}=\omega_{n}^{r}+\omega_{n}^{l}$ with $\omega_{n}^{r}=k$, for some constant $k$, we write 
\begin{align}
M_{ij}&=C_{1}(\epsilon,\lambda_{1})\int_{0}^{2\pi}(k+\omega_{n}^{l})\sin\theta\cos\theta d\theta  \quad i\neq j \\
M_{11}&=C_{2}(\epsilon)\int_{0}^{2\pi}(k+\omega_{n}^{l})\sin^{2}\theta d\theta + C_{3}(\epsilon,\lambda_{1})\int_{0}^{2\pi}(k+\omega_{n}^{l})\cos^{2}\theta d\theta \\
M_{22}&=C_{4}(\epsilon,\lambda_{1})\int_{0}^{2\pi}(k+\omega_{n}^{l})\sin^{2}\theta d\theta + C_{5}(\epsilon)\int_{0}^{2\pi}(k+\omega_{n}^{l})\cos^{2}\theta d\theta,
\end{align}
and $\lambda_{1}^{\prime}(0)$ are determined by the eigenvalues of $M$, with entries given above. Note that the matrix $M$ can also be decomposed as $M=M_{R}+M_{NR}$ where
\begin{align}
M_{R}=\pi(C_{2}+C_{3})kId,
\end{align}
and
\begin{align}
M_{NR}=
\begin{bmatrix}
     C_{2}\int_{0}^{2\pi}\omega_{n}^{l}\sin^{2}\theta d\theta + C_{3}\int_{0}^{2\pi}\omega_{n}^{l}\cos^{2}\theta d\theta   & C_{1}\int_{0}^{2\pi}\omega_{n}^{l}\sin\theta\cos\theta d\theta  \\
    C_{1}\int_{0}^{2\pi}\omega_{n}^{l}\sin\theta\cos\theta d\theta      &  C_{3}\int_{0}^{2\pi}\omega_{n}^{l}\sin^{2}\theta d\theta + C_{2}\int_{0}^{2\pi}\omega_{n}^{l}\cos^{2}\theta d\theta
\end{bmatrix}.
\end{align}
$M_{R}$ being the matrix from the radial term of $V_{n}$ and $M_{NR}$ being the matrix obtained from the length preserving component $\omega_{n}^{l}$. Since $\omega_{n}^{l}$ preserves the length of the inner boundary, we know from \cite{Dambrine} that $B(0,\epsilon)$ provides a critical shape. Furthermore, we have
$$Tr(M_{NR})=0=\lambda_{1,1,NR}^{\prime}(0)+\lambda_{1,2,NR}^{\prime}(0).$$

On the other hand, note that the eigenvalues of $M_{R}$ are given by $k\pi(C_{2}+C_{3})$ with multiplicity two, therefore 
$$\lambda_{1,1,R}^{\prime}(0)=\lambda_{1,2,R}^{\prime}(0)=k\pi(C_{2}+C_{3}).$$
Since $M_{R}$ is a scalar matrix, the eigenvalues of $M=M_{R}+M_{NR}$ are the sum of the eigenvalues of $M_{R}$ and the eigenvalues of $M_{NR}$. Thus, the derivatives of the branches 
$t\mapsto \lambda_{1,1}(t)$, $t\mapsto\lambda_{1,2}(t)$ when we consider $V_{n}=k+\omega_{n}^{l}$ are given by $$\lambda_{1,1}^{\prime}(0)=\lambda_{1,1,NR}^{\prime}(0)+\lambda_{1,1,R}^{\prime}(0) \text{ and }
\lambda_{1,2}^{\prime}(0)=\lambda_{1,2,NR}^{\prime}(0)+\lambda_{1,2,R}^{\prime}(0).$$

Without loss of generality, let $t\mapsto \lambda_{1,1}(t)$ be the smallest of the two branches and let us use the notation $\lambda_{1,1}(t)=\lambda_{1,1}(\Omega_{t})$, with $\Omega=\Omega_{0}$. Noting that $|\partial\Omega_{t}|=2\pi+2\pi(\epsilon-tk)$ and using Lemma \ref{productrule} we have
\begin{align}
\nonumber
\frac{d(\lambda_{1,1}(\Omega_{t})|\partial\Omega_{t}|)}{dt}\Big|_{t=0}=&\frac{d\lambda_{1,1}(\Omega)}{dt}|\partial\Omega|+\lambda_{1,1}(\Omega)\frac{d|\partial\Omega_{t}|}{dt}\Big|_{t=0}\\ \nonumber
=&\lambda_{1,1}^{\prime}(\Omega)|\partial\Omega|+\lambda_{1,1}(\Omega)\frac{d|\partial\Omega_{t}|}{dt}\Big|_{t=0}\\ \nonumber
=& \lambda_{1,1,R}^{\prime}(\Omega)|\partial\Omega| +  \lambda_{1,1,NR}^{\prime}(\Omega)|\partial\Omega| + \lambda_{1,1}(\Omega)\frac{d|\partial\Omega_{t}|}{dt}\Big|_{t=0}\\ \nonumber
=& (\lambda_{1,1,R}^{\prime}(\Omega) +  \lambda_{1,1,NR}^{\prime}(\Omega))2\pi(1+\epsilon) + \lambda_{1,1}(\Omega)\frac{d|\partial\Omega_{t}|}{dt}\Big|_{t=0}.\\ 
\end{align}

Note that if we only deal with $\omega_{n}^{r}=k$, the expression for the Eulerian derivative of the branch of the normalized first Steklov  eigenvalue would be given by
\begin{align}\label{equ581}
\frac{d(\lambda_{1,1}(\Omega_{t})|\partial\Omega_{t}|)}{dt}\Big|_{t=0}=&\frac{d\lambda_{1,1,R}(\Omega)}{dt}2\pi(1+\epsilon)+\lambda_{1,1}(\Omega)(-2\pi k).
\end{align}
Similarly, if we consider only $\omega_{n}^{l}$ then the corresponding expression is given by
\begin{align}\label{equ582}
\frac{d(\lambda_{1,1}(\Omega_{t})|\partial\Omega_{t}|)}{dt}\Big|_{t=0}=&\frac{d\lambda_{1,1,NR}(\Omega)}{dt}2\pi(1+\epsilon).
\end{align}

Let us consider the following function of $\epsilon$.
\begin{align}
E(\epsilon)=&\frac{1+\epsilon^{2}}{2\epsilon(1-\epsilon)}\left(1-\sqrt[]{1-4\epsilon\left(\frac{1-\epsilon}{1+\epsilon^{2}}\right)^{2}}\right)2\pi(1+\epsilon)
\end{align}
Note that $E(\epsilon)\rightarrow 2\pi$ as $\epsilon \rightarrow 0$, $E(\epsilon)\rightarrow 0$ as $\epsilon \rightarrow 1$, and $\epsilon_{1}=0$, $\epsilon_{2}=\frac{-3+\sqrt[]{13}}{2}$ are solutions to $E(\epsilon)=2\pi$. Therefore, since $E(\epsilon)$ is differentiable in $(0,1)$, continuous on $[0,\epsilon_{2}]$ and $E(\epsilon_{1})=E(\epsilon_{2})$, there exists $c_{\epsilon}\in(0,\epsilon_{2})$ such that $\frac{dE(c_{\epsilon})}{d\epsilon}=0$.

On the other hand, note that for the radial deformation (i.e the deformation produced by $\omega_{n}^{r}=k$) we have:
\begin{align}
\nonumber
T_{t}(\vec{x})=&\vec{x}+tV(\vec{x})\\
=&\left\langle \epsilon\cos\theta,\epsilon\sin\theta\right\rangle+t\left\langle -k\cos\theta,-k\sin\theta\right\rangle \\ \nonumber
=&\left\langle \cos\theta,\sin\theta\right\rangle(\epsilon-tk).
\end{align}
This produces  $\Omega_{t}=T_{t}(\Omega)=\mathbb{D}\setminus B(0,\epsilon-tk)$ and we have
\begin{align}
\lambda_{1}(\Omega_{t})=\frac{1}{2(\epsilon-tk)}\frac{1+(\epsilon-tk)^{2}}{1-(\epsilon-tk)}\left(1-\sqrt[]{1-4(\epsilon-tk)\left(\frac{1-(\epsilon-tk)}{1+(\epsilon-tk)^{2}}\right)^{2}}\right).
\end{align}
For the appropriate choice of $k$,  we  have that for the radial variations of the domain
$$\frac{d(\lambda_{1,1}|\partial\Omega_{t}|)}{dt}\Big|_{t=0}=\frac{dE(\epsilon)}{d\epsilon}.$$

That is, the radial variations of the inner domains agree with the variations with respect to $\epsilon$ in the expression for the normalized eigenvalue given by $E(\epsilon)$. This implies by the argument above that there exists an annulus $\mathbb{D}\setminus B(0,c_{\epsilon})$ for which equation \eqref{equ581} satisfies.
\begin{align}\label{equ585}
\frac{d(\lambda_{1,1}|\partial\Omega_{t}|)}{dt}\Big|_{t=0}=0.
\end{align}
In other words, equation \eqref{equ585} implies that $\mathbb{D}\setminus B(0,c_{\epsilon})$ is a critical domain for the case when only radial deformations of the inner boundary are considered.

Based on what we have, consider the admissible class 
$$\CMcal{A}=\{\Omega\subset\mathbb{R}^{2}:\Omega=\mathbb{D}\setminus\Omega_{\delta}, \Omega_{\delta} \text{ and simply connected}\},$$
and also consider the class 
$$\CMcal{B}=\{\Lambda_{\delta}\subset\mathbb{R}^{2}:\Lambda_{\delta}=\mathbb{D}\setminus B(0,\delta)\}.$$
By means of equation \eqref{equ581} we know that there exists $c\in(0,1)$; such that $\Lambda_{c}$ provides a critical domain for $\lambda_{1,1}(\Lambda_{\delta})|\partial\Lambda_{\delta}|$ among sets belonging to $\mathcal{B}$. Furthermore, for those $\Lambda_{\delta}$ with $|\partial\Lambda_{\delta}|=2\pi c$, equation \eqref{equ582} guarantees that among the sets $\mathbb{D}\setminus \Omega_{\delta}$ in  $\mathcal{A}$ with $|\partial\Omega_{\delta}|=2\pi c$; $\Lambda_{c}$ is also a critical shape for the shape functional $\lambda_{1,1}(\Omega)|\partial\Omega|$. We therefore have proved the following result.
\begin{thm}\label{mainresult}
For each $\epsilon\in(0,1)$, $\Lambda_{\epsilon}=\mathbb{D}\setminus B(0,\epsilon)$ provides a critical domain among the class of annular domains $\Omega=\mathbb{D}\setminus \Omega_{\delta}$, with $\Omega_{\delta}$ simply connected subset of $\mathbb{D}$ and $|\partial\Omega_{\delta}|=2\pi\epsilon$.
\end{thm}

\textbf{Remark:}
Note that in proving Theomre \ref{mainresult}, we have used the following fact. If $f,h,g$ are differentiable functions such $f=h+g$, and $h^{\prime}(a)=0$ and $g^{\prime}(a)=0$, then $f^{\prime}(a)=0$. Noting that $\Lambda_{c}$ provides a critical domain from equation \eqref{equ581}, and it also provides a critical domain from equation \eqref{equ582}, and since the sum equations \eqref{equ581} and \eqref{equ582} gives 
\begin{align}
\frac{d(\lambda_{1,1}(\Omega_{t})|\partial\Omega_{t}|)}{dt}.
\end{align}
Then, $\Lambda_{c}$ is a critical domain for the shape functional $\lambda_{1,1}(\Omega)|\partial\Omega|$. In other words we have

\begin{corr}\label{propoepsilon}
The annulus $\Lambda_{\epsilon_{0}}=\mathbb{D}\setminus B(0,\epsilon_{0})$, where $\epsilon_{0}$ is the unique root in $(0,1)$ of the polynomial $\Pi(\epsilon)=\epsilon^{6}-10\epsilon^{5}+23\epsilon^{4}-12\epsilon^{3}+23\epsilon^{2}-10\epsilon+1$ is a critical domain in $\CMcal{A}$.
\end{corr}

\section{Numerical Experiments and Discussion}
In this section we present a series of numerical results obtained with FreeFem++ \cite{freefem}. We used a slight variation of the code created by Bogosel and presented in \cite{Bogosel}. These examples allow us to verify the results obtained in the previous section, and to conjecture that indeed an annulus of a given perimeter maximizes the first normalized Steklov eigenvalue. In particular, that the annulus of inner radius given by $\epsilon_{0}$ in Corollary \ref{propoepsilon} is a candidate for a global maximizer (up to scalings of the domains). The conjecture is the following. 
\begin{conj}[\cite{Girouard3}]
When restricting to bounded connected planar domains with two boundary components, the expectation is that the best planar annulus is the one that realizes the max on the curve of the function $\lambda_{1}(\Lambda_{\epsilon})|\partial\Lambda_{\epsilon}|$, where $\lambda_{1}(\Lambda_{\epsilon})$ is given by equation (\ref{steklovevals}).
\end{conj}

\subsection{Translations of the inner circle along the $x$-axis}
In these examples we translate the center of inner circle through the $x$-axis, we compute the first normalized eigenvalue for different values of the inner radius. Tables \ref{table1}, \ref{table2} and \ref{table3} show that the maximum values are attained when the inner circle is centered at the origin. In particular, the maximum value for these examples is attained at the annulus with inner radius $\epsilon_{0}$ described by Corollary \ref{propoepsilon}.
\begin{table}[H]
\begin{center}
     \raisebox{-0.5\totalheight}{\includegraphics[width=0.4\textwidth]{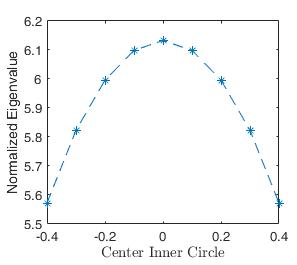}}
    \begin{tabular}{|c|c|c}
    \hline
      Center   &  $2\pi(1+\epsilon)\lambda_{1}(\Lambda_{\epsilon})$ \\
      \hline
      $(-0.4,0)$   &  $5.5724$\\
      \hline
      $(-0.3,0)$  &   $5.8231$\\
      \hline
      $(-0.2,0)$  &   $5.9960$\\
      \hline
      $(-0.1,0)$  &   $6.0987$\\
      \hline
      $\textbf{(0,0)}$     &   $\textbf{6.1328}$\\
      \hline
      $\ \ (0.1,0)$  &   $6.0987$\\
      \hline
      $\ \ (0.2,0)$  &   $5.9960$\\
      \hline
      $\ \ (0.3,0)$  &   $5.8231$\\
      \hline
      $\ \ (0.4,0)$  &   $5.5724$\\
      \hline
      \end{tabular} 
      \caption{First normalized eigenvalues $\lambda_{1}(\Lambda_{\epsilon})$, with $\epsilon=0.3$}
      \label{table1}
      \end{center}
\end{table}

\begin{table}[H]
\begin{center}
     \raisebox{-0.5\totalheight}{\includegraphics[width=0.4\textwidth]{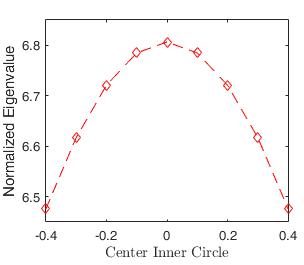}}
    \begin{tabular}{|c|c|c}
    \hline
      Center   &  $2\pi(1+\epsilon)\lambda_{1}(\Lambda_{\epsilon})$ \\
      \hline
      $(-0.4,0)$   &  $6.4759$\\
      \hline
      $(-0.3,0)$  &   $6.6169$\\
      \hline
      $(-0.2,0)$  &   $6.7208$\\
      \hline
      $(-0.1,0)$  &   $6.7848$\\
      \hline
      $\textbf{(0,0)}$     &   $\textbf{6.8064}$\\
      \hline
      $\ \ (0.1,0)$  &   $6.7848$\\
      \hline
      $\ \ (0.2,0)$  &   $6.7208$\\
      \hline
      $\ \ (0.3,0)$  &   $6.6169$\\
      \hline
      $\ \ (0.4,0)$  &   $6.4759$\\
      \hline
      \end{tabular} 
      \caption{First normalized eigenvalues $\lambda_{1}(\Lambda_{\epsilon})$, with $\epsilon=0.146721$}
      \label{table2}
      \end{center}
\end{table}

\begin{table}[H]
\begin{center}
     \raisebox{-0.5\totalheight}{\includegraphics[width=0.4\textwidth]{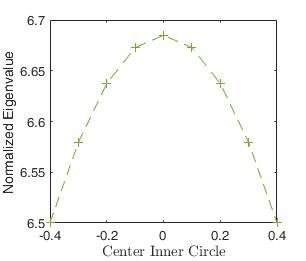}}
    \begin{tabular}{|c|c|c}
    \hline
      Center   &  $2\pi(1+\epsilon)\lambda_{1}(\Lambda_{\epsilon})$ \\
      \hline
      $(-0.4,0)$   &  $6.5001$\\
      \hline
      $(-0.3,0)$  &   $6.5794$\\
      \hline
      $(-0.2,0)$  &   $6.6374$\\
      \hline
      $(-0.1,0)$  &   $6.6729$\\
      \hline
      $\textbf{(0,0)}$     &   $\textbf{6.6849}$\\
      \hline
      $\ \ (0.1,0)$  &   $6.6729$\\
      \hline
      $\ \ (0.2,0)$  &   $6.6374$\\
      \hline
      $\ \ (0.3,0)$  &   $6.5794$\\
      \hline
      $\ \ (0.4,0)$  &   $6.5001$\\
      \hline
      \end{tabular} 
      \caption{First normalized eigenvalues $\lambda_{1}(\Lambda_{\epsilon})$, with $\epsilon=0.08$}
      \label{table3}
      \end{center}
\end{table}
\subsection{Translations of the inner circle along the line $y=-x$}
We also considered examples in which the center of the inner circle moves along the line $y=-x$, we considered the same values for $\epsilon$ as in the case of the center moving along the $x$-axis. We can see that moving the center of the inner circle along $y=-x$  produces slightly bigger values than the horizontal movement along the $x-$axis, but at the origin the values are the same in both cases. Here again, we obtain a maximum value for the annulus with inner radius $\epsilon_{0}$ determined in Corollary \ref{propoepsilon}. Tables, \ref{table4} \ref{table5} and \ref{table6} show the corresponding numerical values.
\begin{table}[H]
\begin{center}
     \raisebox{-0.5\totalheight}{\includegraphics[width=0.4\textwidth]{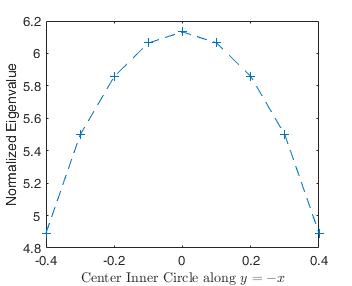}}
    \begin{tabular}{|c|c|c}
    \hline
      Center   &  $2\pi(1+\epsilon)\lambda_{1}(\Lambda_{\epsilon})$ \\
      \hline
      $(-0.4,0.4)$   &  $4.8916$\\
      \hline
      $(-0.3,0.3)$  &   $5.4976$\\
      \hline
      $(-0.2,0.2)$  &   $5.8580$\\
      \hline
      $(-0.1,0.1)$  &   $6.0645$\\
      \hline
      $\textbf{(0,0)}$     &   $\textbf{6.1328}$\\
      \hline
      $\ \ (0.1,-0.1)$  &   $6.0645$\\
      \hline
      $\ \ (0.2,-0.2)$  &   $5.8580$\\
      \hline
      $\ \ (0.3,-0.3)$  &   $5.4976$\\
      \hline
      $\ \ (0.4,-0.4)$  &   $4.8916$\\
      \hline
      \end{tabular} 
      \caption{First normalized eigenvalues $\lambda_{1}(\Lambda_{\epsilon})$, with $\epsilon=0.3$}
      \label{table4}
      \end{center}
\end{table}

\begin{table}[H]
\begin{center}
     \raisebox{-0.5\totalheight}{\includegraphics[width=0.4\textwidth]{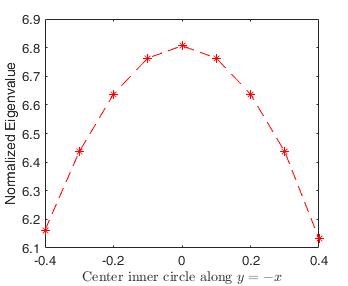}}
    \begin{tabular}{|c|c|c}
    \hline
      Center   &  $2\pi(1+\epsilon)\lambda_{1}(\Lambda_{\epsilon})$ \\
      \hline
      $(-0.4,0.4)$   &  $6.1623$\\
      \hline
      $(-0.3,0.3)$  &   $6.4363$\\
      \hline
      $(-0.2,0.2)$  &   $6.6375$\\
      \hline
      $(-0.1,0.1)$  &   $6.7633$\\
      \hline
      $\textbf{(0,0)}$     &   $\textbf{6.8064}$\\
      \hline
      $\ \ (0.1,-0.1)$  &   $6.7633$\\
      \hline
      $\ \ (0.2,-0.2)$  &   $6.6375$\\
      \hline
      $\ \ (0.3,-0.3)$  &   $6.4363$\\
      \hline
      $\ \ (0.4,-0.4)$  &   $6.1623$\\
      \hline
      \end{tabular} 
      \caption{First normalized eigenvalues $\lambda_{1}(\Lambda_{\epsilon})$, with $\epsilon=0.146721$}
      \label{table5}
      \end{center}
\end{table}

\begin{table}[H]
\begin{center}
     \raisebox{-0.5\totalheight}{\includegraphics[width=0.4\textwidth]{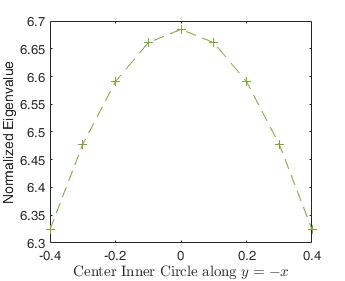}}
    \begin{tabular}{|c|c|c}
    \hline
      Center   &  $2\pi(1+\epsilon)\lambda_{1}(\Lambda_{\epsilon})$ \\
      \hline
      $(-0.4,0.4)$   &  $6.3244$\\
      \hline
      $(-0.3,0.3)$  &   $6.4777$\\
      \hline
      $(-0.2,0.2)$  &   $6.5909$\\
      \hline
      $(-0.1,0.1)$  &   $6.6610$\\
      \hline
      $\textbf{(0,0)}$     &   $\textbf{6.6849}$\\
      \hline
      $\ \ (0.1,-0.1)$  &   $6.6610$\\
      \hline
      $\ \ (0.2,-0.2)$  &   $6.5909$\\
      \hline
      $\ \ (0.3,-0.3)$  &   $6.4777$\\
      \hline
      $\ \ (0.4,-0.4)$  &   $6.3244$\\
      \hline
      \end{tabular} 
      \caption{First normalized eigenvalues $\lambda_{1}(\Lambda_{\epsilon})$, with $\epsilon=0.08$}
      \label{table6}
      \end{center}
\end{table}

\subsection{Radial perturbation of the inner boundary}
For these examples, we use the fact that the length $L$ of a curve in polar coordinates is given by $L=\int_{\alpha}^{\beta}\sqrt{r^{2}+(dr/d\theta)^{2}}d\theta$. By considering the function $r=acos(k\theta)+b$ with $\alpha=0$ and $\beta=2\pi$, we can find values of the parameters $a$,$b$ and $k$ (with $k$ being an integer) such that $L=2\pi\epsilon$. 

Solving the corresponding integral with the given function, we find that 
\begin{align}\label{length1}
     L=a^{2}\pi+2\pi b^{2}+a^{2}k^{2}\pi.  
\end{align}
Then, $L=2\pi\epsilon$ implies $\epsilon = \frac{a^{2}}{2}+b^{2}+\frac{a^{2}k^{2}}{2}$. Since $b$ determines the maximum and minimum values for $r$. We choose $b=\epsilon_{0}$, so (\ref{length1}) implies that we must choose $\epsilon\geq\epsilon_{0}^{2}\approx 0.021527$. In particular, if we allow $\epsilon = \epsilon_{0}$, again from (\ref{length1}) we obtain the corresponding expression for the parameter $a$ given by
\begin{align}
a = \left|\sqrt{\frac{2(\epsilon_{0}-\epsilon_{0}^{2})}{1+k^{2}}} \right|.
\end{align}
With these values at hand we use FreeFem++ to obtain the corresponding numerical values for the normalized eigenvalues. Table \ref{table7} show different values for different domains, all of which have perimeter $2\pi\epsilon_{0}$.
\begin{table}[H]
\begin{center}
    \begin{tabular}{|c|c|c|}
    \hline
      k   &  $a\approx$ & $2\pi(1+\epsilon)\lambda_{1}(\Omega)$ \\
      \hline
      $5$ & $0.0981$ & $6.0338$ \\
      \hline
      $10$ & $0.0497$ & $6.3146$ \\
      \hline 
      $20$ & $0.0249$ & $6.4700$ \\
      \hline
      $50$ & $0.01$ & $6.5698$ \\
      \hline
      \end{tabular} 
      \caption{First normalized eingenvalues, for different domains $\Omega$ satisfying  $|\partial\Omega|= 2\pi\lambda_{1}(\Omega)(1+\epsilon_{0})$.}
      \label{table7}
      \end{center}
\end{table}
Notice that the numerical values showed in Table \ref{table7} are all less or equal than the value $6.8064$ obtained when the inner domain is the circle of radius $\epsilon_{0}$. It is worth to mention that one can choose different values for $k$ with the corresponding value of $a$ associated to it. When choosing the values for $k$, one have to be careful so that FreeFem++ does not produces errors of crossed boundaries. Examples of domains with $k=5, 10, 20$ and $50$ are shown in Figure 2. 
\begin{table}[H]
\begin{center}
    \begin{tabular}{c c}
      \raisebox{-0.5\totalheight}{\includegraphics[width=0.4\textwidth]{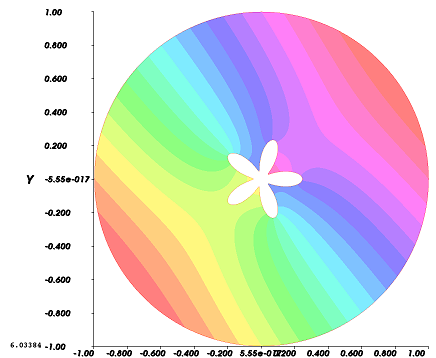}} & \raisebox{-0.5\totalheight}{\includegraphics[width=0.4\textwidth]{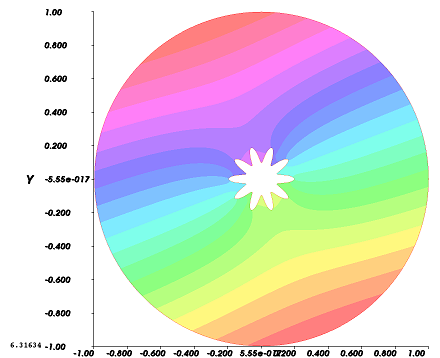}} \\ 
      \  \ & \ \ \\ 
      \raisebox{-0.5\totalheight}{\includegraphics[width=0.4\textwidth]{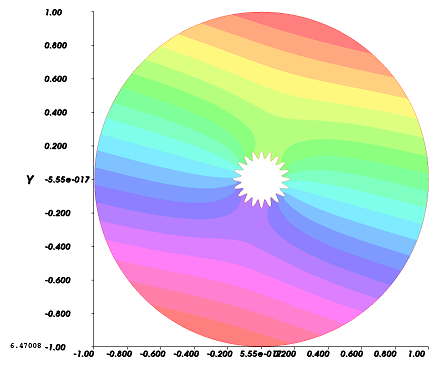}} & \raisebox{-0.5\totalheight}{\includegraphics[width=0.4\textwidth]{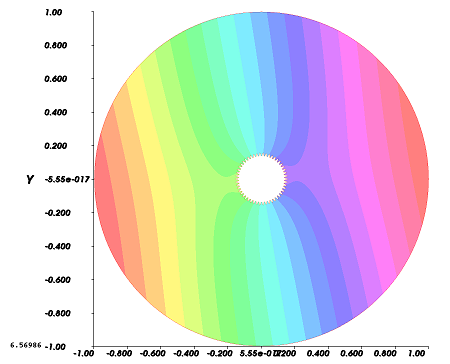}}
      \end{tabular} 
      \caption*{Figure 2: Different domains satisfying $|\partial\Omega|=2\pi\lambda_{1}(\Omega)(1+\epsilon_{0})$}
      \label{table8}
      \end{center}
\end{table}

\subsection{Discussion}
Although more numerical examples can be cooked up to test for maximality of the proposed domain; a proof that the domain described in Corollary \ref{propoepsilon} is the maximizer for the first normalized Steklov eigenvalue, among planar annular domains of fixed perimeter is still missing. The results obtained here are local in the sense that they are obtained only in the direction of the vector field $V$ from Theorem \ref{thm1}. Attempts to compute the second variation of the shape functional in the direction of $V$ were made. However, the author was not able to get any definitive conclusion following this approach. Perhaps a more analytical approach will provide a global result.


  \textbf{Aknowledgements:} The author would like to acknowledge  the support by Pontificia Universidad Javeriana at Bogotá, D.C under the project ID-000000000010592. In addition, the author would like to thank Dr. Alexandre Girouard for pointing out the question that have origin to this manuscript.

\bibliography{REFS}
}\end{document}